\documentclass[12pt]{article}
\usepackage{amsmath}
\usepackage{amsfonts}
\usepackage{amssymb}

\def\thebibliograph#1#2{\section*{{\normalsize \bf #2}}\list
   {[\arabic{enumi}]}{\settowidth\labelwidth{[#1]}\leftmargin\labelwidth
     \advance\leftmargin\labelsep
     \usecounter{enumi}}
     \def\newblock{\hskip .11em plus .33em minus -.07em}
     \sloppy
     \sfcode`\.=1000\relax}

\newtheorem{theorem}{Theorem}

\newtheorem{lemma}{Lemma}

\input{tcilatex}
\begin{document}

\title{A note on Hardy-type inequalities in variable exponent Lebesgue spaces%
}
\author{Douadi Drihem}
\date{\today }
\maketitle

\begin{abstract}
We present new estimate for Hardy-type\ inequality in variable exponent
Lebesgue spaces. More precisely, by imposing regularity assumptions on the
exponent, we prove that the estimations can be reduced to the fixed
exponents.


\end{abstract}

\section{Introduction}

It is well known that Hardy-type inequalities\ play an important role in
Harmonic Analysis. For instance, they appear in the interpolation of spaces\ 
$\mathrm{\cite{BS88}}$, in the study of hyperbolic partial differential
equations and for studying the decay of linear waves on black hole back $%
\mathrm{\cite{DF10}}$. We can find some interesting applications of
Hardy-type inequalities in $\mathrm{\cite{KP03}}$ and references therein.

The classical Hardy inequalities says that%
\begin{equation*}
\Big\|t^{s}\int_{t}^{\infty }\tau ^{-s}\varepsilon _{\tau }\frac{d\tau }{%
\tau }\Big\|_{L^{q}((0,\infty ),\frac{dt}{t})}+\Big\|t^{-s}\int_{0}^{t}\tau
^{s}\varepsilon _{\tau }\frac{d\tau }{\tau }\Big\|_{L^{q}((0,\infty ),\frac{%
dt}{t})}\lesssim \big\|\varepsilon _{\tau }\big\|_{L^{q}((0,\infty ),\frac{dt%
}{t})}
\end{equation*}%
for any $s>0$ and any $1\leq q<\infty $. This statement in variable exponent
Lebesgue spaces was first proved by V. Kokilashvili and S. Samko \cite%
{KSamko04}\ and by L. Diening and S. Samko \cite{D.S2007} under the
assumption that $q$\textit{\ }is \emph{$\log $}-H\"{o}lder continuous both
at the origin\ and at infinity. More results for Hardy-type inequalities in
variable exponent Lebesgue spaces can be found in S. Boza and J. Soria $%
\mathrm{\cite{Bs07}}$,\textrm{\ }Cruz-Uribe and Mamedov $\mathrm{\cite{CM12}}
$,\textrm{\ }P. Harjulehto, P. H\"{a}st\"{o}, and M. Koskinoja $\mathrm{\cite%
{HHK05}}$,\textrm{\ }F. I. Mamedov and A. Harman $\mathrm{\cite{M10}}$, and
references therein. Here under the same assumptions we prove that $\big\|%
\varepsilon _{\tau }\big\|_{L^{q(\cdot )}((0,\infty ),\frac{dt}{t})}$ can be
replaced by 
\begin{equation*}
\Big(\int_{0}^{1}\varepsilon _{t}^{q(0)}\frac{dt}{t}\Big)^{\frac{1}{q(0)}}+%
\Big(\int_{1}^{\infty }\varepsilon _{t}^{q_{\infty }}\frac{dt}{t}\Big)^{%
\frac{1}{q_{\infty }}}.
\end{equation*}

\section{Preliminaries}

As usual, we denote by $\mathbb{R}$ the reals, $\mathbb{N}$ the collection
of all natural numbers and $\mathbb{N}_{0}=\mathbb{N}\cup \{0\}$. The letter 
$\mathbb{Z}$ stands for the set of all integer numbers.\ The expression $%
f\lesssim g$ means that $f\leq c\,g$ for some independent constant $c$ (and
non-negative functions $f$ and $g$), and $f\approx g$ means $f\lesssim
g\lesssim f$.\vskip5pt

If $E\subset {\mathbb{R}}$ is a measurable set, then $\chi _{E}$ denotes the
characteristic function of the set $E$.\vskip5pt

By $c$ we denote generic positive constants, which may have different values
at different occurrences. Although the exact values of the constants are
usually irrelevant for our purposes, sometimes we emphasize their dependence
on certain parameters (e.g. $c(p)$ means that $c$ depends on $p$, etc.).
Further notation will be properly introduced whenever needed.

The variable exponents that we consider are always measurable functions $p$
on $\mathbb{R}$ with range in $[1,\infty \lbrack $. We denote the set of
such functions by $\mathcal{P}({\mathbb{R}})$. We use the standard notation $%
p^{-}:=\underset{x\in \mathbb{R}}{\text{ess-inf}}$ $p(x)$,$\quad p^{+}:=%
\underset{x\in \mathbb{R}}{\text{ess-sup }}p(x)$.

The variable exponent modular is defined by 
\begin{equation*}
\varrho _{p(\cdot )}(f):=\int_{\mathbb{R}}\varrho _{p(x)}(\left\vert
f(x)\right\vert )dx,
\end{equation*}%
where $\varrho _{p}(t)=t^{p}$. The variable exponent Lebesgue space $%
L^{p(\cdot )}$\ consists of measurable functions $f$ on $\mathbb{R}$ such
that $\varrho _{p(\cdot )}(\lambda f)<\infty $ for some $\lambda >0$. We
define the Luxemburg norm on this space by the formula 
\begin{equation*}
\left\Vert f\right\Vert _{p(\cdot )}:=\inf \Big\{\lambda >0:\varrho
_{p(\cdot )}\Big(\frac{f}{\lambda }\Big)\leq 1\Big\}.
\end{equation*}%
A useful property is that $\left\Vert f\right\Vert _{p(\cdot )}\leq 1$ if
and only if $\varrho _{p(\cdot )}(f)\leq 1$, see \cite{DHHR}, Lemma 3.2.4.

Let $p,q\in \mathcal{P}(\mathbb{R})$. The mixed Lebesgue-sequence space $%
\ell ^{q(\cdot )}(L^{p(\cdot )})$ is defined on sequences of $L^{p(\cdot )}$%
-functions by the modular%
\begin{equation*}
\varrho _{\ell ^{q(\cdot )}(L^{p\left( \cdot \right)
})}((f_{v})_{v}):=\sum\limits_{v=-\infty }^{\infty }\inf \Big\{\lambda
_{v}>0:\varrho _{p(\cdot )}\Big(\frac{f_{v}}{\lambda _{v}^{1/q(\cdot )}}\Big)%
\leq 1\Big\}.
\end{equation*}%
The (quasi)-norm is defined from this as usual:%
\begin{equation}
\left\Vert \left( f_{v}\right) _{v}\right\Vert _{\ell ^{q(\cdot
)}(L^{p\left( \cdot \right) })}:=\inf \Big\{\mu >0:\varrho _{\ell ^{q(\cdot
)}(L^{p(\cdot )})}\Big(\frac{1}{\mu }(f_{v})_{v}\Big)\leq 1\Big\}.
\label{mixed-norm}
\end{equation}%
If $q^{+}<\infty $, then we can replace $\mathrm{\eqref{mixed-norm}}$ by the
simpler expression $\varrho _{\ell ^{q(\cdot )}(L^{p(\cdot
)})}((f_{v})_{v}):=\sum\limits_{v=-\infty }^{\infty }\left\Vert
|f_{v}|^{q(\cdot )}\right\Vert _{\frac{p(\cdot )}{q(\cdot )}}$.

We say that a function $g\,:\,{\mathbb{R}}\rightarrow \mathbb{R}$ is \emph{$%
\log $-H\"{o}lder continuous at the origin} (or has a \emph{$\log $ decay at
the origin}), if there exists a constant $c_{\log }(g)>0$ such that 
\begin{equation*}
|g(x)-g(0)|\leq \frac{c_{\log }(g)}{\ln (e+1/|x|)}
\end{equation*}%
for all $x\in \mathbb{R}$. If, for some $g_{\infty }\in \mathbb{R}$ and $%
c_{\log }>0$, there holds 
\begin{equation*}
|g(x)-g_{\infty }|\leq \frac{c_{\log }}{\ln (e+|x|)}
\end{equation*}%
for all $x\in \mathbb{R}$, then we say that $g$ is \emph{$\log $-H\"{o}lder
continuous at infinity} (or has a \textit{$\log $ decay at infinity}).%
\textit{\ }The constants $c_{\log }(g)$ and $c_{\log }$ are called the 
\textit{locally }log\textit{-H\"{o}lder constant }and the log\textit{-H\"{o}%
lder decay constant}, respectively. We refer to the recent monograph $%
\mathrm{\cite{CF13}}$ for further properties, historical remarks and
references on variable exponent spaces.

\subsection{Technical lemmas}

In this subsection we present some results which are useful for us. The
following lemma is a Hardy-type inequality which is easy to prove.

\begin{lemma}
\label{Hardy-inequality}\textit{Let }$0<a<1$, $\sigma \geq 0$ \textit{and }$%
0<p\leq \infty $\textit{. Let }$\left\{ \varepsilon _{k}\right\} _{k}$%
\textit{\ be a sequences of positive real numbers and denote} $\delta
_{k}=\sum_{j=-\infty }^{\infty }\left\vert k-j\right\vert ^{\sigma
}a^{\left\vert k-j\right\vert }\varepsilon _{j}$.\textit{\ }Then there
exists constant $c>0\ $\textit{depending only on }$a$\textit{\ and }$p$ such
that%
\begin{equation*}
\Big(\sum\limits_{k=-\infty }^{\infty }\delta _{k}^{p}\Big)^{1/p}\leq c\text{
}\Big(\sum\limits_{k=-\infty }^{\infty }\varepsilon _{k}^{p}\Big)^{1/p}.
\end{equation*}
\end{lemma}

We will make use of the following statement, see \cite{DHHMS}, Lemma 3.3 for 
$w:=1$.

\begin{lemma}
\label{DHHR-estimate}Let $Q=(a,b)\subset \mathbb{R}$ with $0<a<b<\infty $.
Let $p\in \mathcal{P}(\mathbb{R})$ \textit{be log-H\"{o}lder continuous at
the origin }and $w:\mathbb{R}\rightarrow \mathbb{R}^{+}$ be a weight
function. Then for every $m>0$ there exists $\gamma =e^{-4mc_{\log
}(1/p)}\in \left( 0,1\right) $ such that%
\begin{eqnarray*}
&&\Big(\frac{\gamma }{w(Q)}\int_{Q}\left\vert f(y)\right\vert w(y)dy\Big)%
^{p\left( x\right) } \\
&\leq &\max \big(1,\left( w(Q)\right) ^{1-\frac{p\left( x\right) }{p^{-}}}%
\big)\frac{1}{w(Q)}\int_{Q}\left\vert f(y)\right\vert ^{p\left( y,0\right)
}w(y)dy \\
&&+\omega (m,b)\Big(\frac{1}{w(Q)}\int_{Q}g(x,y)w(y)dy\Big)
\end{eqnarray*}%
hold if $0<w(Q)<\infty $, all $x\in Q\subset \mathbb{R}$ and all $f\in
L^{p\left( \cdot \right) }(w)+L^{\infty }$\ with $\left\Vert fw^{1/p\left(
\cdot \right) }\right\Vert _{p\left( \cdot \right) }+\left\Vert f\right\Vert
_{\infty }\leq 1$, where%
\begin{equation*}
\omega (m,b)=\min \left( b^{m},1\right) \text{, }p\left( y,0\right) =p\left(
y\right) \text{ \ \ and \ \ }g(x,y)=(e+\frac{1}{x})^{-m}+(e+\frac{1}{y})^{-m}
\end{equation*}%
or%
\begin{equation*}
\omega (m,b)=\min \left( b^{m},1\right) \text{, }p\left( y,0\right) =p\left(
0\right) \text{ \ \ and \ \ }g(x,y)=(e+\frac{1}{x})^{-m}\chi _{\{x:p(x){<}%
p(0)\}}(x).
\end{equation*}%
\textit{In addition we have the same estimate}, where%
\begin{equation*}
\omega (m,b)=1\text{, }p\left( y,0\right) =p_{\infty }\text{ \ \ and \ \ }%
g(x,y)=(e+x)^{-m}\chi _{\{x:p(x){<}p_{\infty }\}}(x),
\end{equation*}%
if $p\in \mathcal{P}(\mathbb{R})$ satisfies the log\textit{-H\"{o}lder decay
condition, where we take }$\gamma =e^{-4mc_{\log }}$.
\end{lemma}

Notice that in the proof of this theorem we need only that 
\begin{equation*}
\int_{Q}\left\vert f(y)\right\vert ^{p\left( y\right) }w(y)dy\leq 1
\end{equation*}%
and/or $\left\Vert f\right\Vert _{\infty }\leq 1$. The proof of this lemma\
is given in \cite{D3}.

\section{Main results}

Various important results have been proved in the space $L^{p(\cdot )}$
under some assumptions on $p$ such us the boundedness of the maximal
operator in $L^{p(\cdot )}$ spaces on bounded domains. This fact was first
realized by L. Diening \cite{Di}. This statement was then extended to the
unbounded case by D. Cruz-Uribe, A. Fiorenza and C. Neugebauer \cite{CUFN03}%
. Estimates for potential type operators in variable $L^{p(\cdot )}$ spaces
were first considered by Samko \cite{Sam98}. Fractional maximal operators
were first studied in this setting by Kokilashvili and Samko \cite{KSam03}.
We refer to \cite{DHHR} for further contributions and historical remarks in
the study of singular integral and fractional integral operators in variable
exponent spaces. As mentioned in the introduction\ we present some new
estimate for Hardy operators $\int_{t}^{\infty }\tau ^{-s}\varepsilon _{\tau
}\frac{d\tau }{\tau }$ and $\int_{0}^{t}\tau ^{s}\varepsilon _{\tau }\frac{%
d\tau }{\tau }$, $t>0$\textit{\ }in variable exponent Lebesgue spaces. More
precisely, we have the following results:

\begin{theorem}
\label{main-result}\textit{Let }$s>0$\textit{. Let }$p\in \mathcal{P}(%
\mathbb{R})$\textit{\ }be \emph{$\log $}-H\"{o}lder continuous both at the
origin\ and at infinity with $1\leq p^{-}\leq p^{+}<\infty $\textit{. Let }$%
\left\{ \varepsilon _{t}\right\} _{t}$\textit{\ be a sequence of positive
measurable functions.} Let%
\begin{equation*}
\eta _{t}=t^{s}\int_{t}^{\infty }\tau ^{-s}\varepsilon _{\tau }\frac{d\tau }{%
\tau }\quad \text{and\quad }\lambda _{t}=t^{-s}\int_{0}^{t}\tau
^{s}\varepsilon _{\tau }\frac{d\tau }{\tau }.
\end{equation*}%
Then there exists constant $c>0\ $\textit{depending\ only\ on\ }$s$, $p^{-}$%
\textit{,\ c}$_{\log }(p)$ \textit{and }$p^{+}$\ such that%
\begin{equation}
\left\Vert \eta _{t}\right\Vert _{L^{p(\cdot )}((0,\infty ),\frac{dt}{t}%
)}\approx \Big(\int_{0}^{1}\eta _{t}^{p(0)}\frac{dt}{t}\Big)^{\frac{1}{p(0)}%
}+\Big(\int_{1}^{\infty }\eta _{t}^{p_{\infty }}\frac{dt}{t}\Big)^{\frac{1}{%
p_{\infty }}}  \label{Est1.1.1}
\end{equation}%
and%
\begin{equation}
\left\Vert \lambda _{t}\right\Vert _{L^{p(\cdot )}((0,\infty ),\frac{dt}{t}%
)}\approx \Big(\int_{0}^{1}\lambda _{t}^{p(0)}\frac{dt}{t}\Big)^{\frac{1}{%
p(0)}}+\Big(\int_{1}^{\infty }\lambda _{t}^{p_{\infty }}\frac{dt}{t}\Big)^{%
\frac{1}{p_{\infty }}}.  \label{Est1.1.2}
\end{equation}%
Moreover, 
\begin{equation*}
\left\Vert \eta _{t}\right\Vert _{L^{p(\cdot )}((0,\infty ),\frac{dt}{t}%
)}+\left\Vert \lambda _{t}\right\Vert _{L^{p(\cdot )}((0,\infty ),\frac{dt}{t%
})}\lesssim \Big(\int_{0}^{1}\varepsilon _{t}^{p(0)}\frac{dt}{t}\Big)^{\frac{%
1}{p(0)}}+\Big(\int_{1}^{\infty }\varepsilon _{t}^{p_{\infty }}\frac{dt}{t}%
\Big)^{\frac{1}{p_{\infty }}}.
\end{equation*}
\end{theorem}

\textbf{Proof. }We will do the proof in several steps.

\textit{Step 1.} We prove that%
\begin{equation*}
\left\Vert \eta _{t}\right\Vert _{L^{p(\cdot )}((0,\infty ),\frac{dt}{t}%
)}\lesssim \Big(\int_{0}^{1}\eta _{t}^{p(0)}\frac{dt}{t}\Big)^{\frac{1}{p(0)}%
}+\Big(\int_{1}^{\infty }\eta _{t}^{p_{\infty }}\frac{dt}{t}\Big)^{\frac{1}{%
p_{\infty }}}.
\end{equation*}%
We suppose that the right-hand side is less than or equal one. Notice that%
\begin{equation*}
\left\Vert \eta _{t}\right\Vert _{L^{p(\cdot )}((0,\infty ),\frac{dt}{t}%
)}\approx \Big\|\Big(t^{-\frac{1}{p(t)}}\eta _{t}\chi _{\lbrack
2^{v},2^{1+v}]}\Big)_{v}\Big\|_{\ell ^{p(\cdot )}(L^{p(\cdot )})}.
\end{equation*}%
We see that%
\begin{equation*}
\eta _{t}=\frac{t^{s}}{\log 2}\int_{\frac{t}{2}}^{t}t^{-s}\eta _{t}\frac{%
d\tau }{\tau }\lesssim \frac{t^{s}}{\log 2}\int_{\frac{t}{2}}^{t}\tau
^{-s}\eta _{\tau }\frac{d\tau }{\tau }\leq \frac{t^{s}}{\log 2}%
\int_{2^{v-1}}^{\infty }\tau ^{-s}\eta _{\tau }\frac{d\tau }{\tau }
\end{equation*}%
for any $v\in \mathbb{Z}$ and any $t\in \lbrack 2^{v},2^{v+1}]$. We write, 
\begin{equation*}
t^{s}\int_{2^{v-1}}^{\infty }\tau ^{-s}\eta _{\tau }\frac{d\tau }{\tau }%
=\sum_{j=v-1}^{\infty }t^{s}\int_{2^{j}}^{2^{j+1}}\tau ^{-s}\eta _{\tau }%
\frac{d\tau }{\tau },\quad v\in \mathbb{Z}.
\end{equation*}%
\textit{Substep 1.1.} $v\leq 0$. For any $t\in \lbrack 2^{v},2^{v+1}]$, we
write%
\begin{equation*}
\eta _{t,v}=\sum_{j=v-1}^{\infty }t^{s}\int_{2^{j}}^{2^{j+1}}\tau ^{-s}\eta
_{\tau }\frac{d\tau }{\tau }=\sum_{j=v-1}^{-1}\cdot \cdot \cdot
+\sum_{j=0}^{\infty }\cdot \cdot \cdot =\eta _{t,1,v}+\eta _{t,2,v}.
\end{equation*}%
\textit{Estimation of }$\eta _{t,1,v}$. Let $\sigma >0$ be such that $%
p^{+}<\sigma $.\ We have%
\begin{eqnarray*}
\Big(\sum_{j=v-1}^{-1}\int_{2^{j}}^{2^{j+1}}\tau ^{-s}\eta _{\tau }\frac{%
d\tau }{\tau }\Big)^{p(t)/\sigma } &\leq &\sum_{j=v-1}^{-1}\Big(%
\int_{2^{j}}^{2^{j+1}}\tau ^{-s}\eta _{\tau }\frac{d\tau }{\tau }\Big)%
^{p(t)/\sigma } \\
&\leq &\sum_{j=v-1}^{-1}2^{-jsp(t)/\sigma }\Big(\int_{2^{j}}^{2^{j+1}}\eta
_{\tau }\frac{d\tau }{\tau }\Big)^{p(t)/\sigma } \\
&=&2^{-vsp(t)/\sigma }\sum_{j=v-1}^{-1}2^{(v-j)sp(t)/\sigma }\Big(%
\int_{2^{j}}^{2^{j+1}}\eta _{\tau }\frac{d\tau }{\tau }\Big)^{p(t)/\sigma }.
\end{eqnarray*}%
By H\"{o}lder's inequality, we estimate this expression by%
\begin{equation*}
c2^{-vsp(t)/\sigma }\Big(\sum_{j=v-1}^{-1}2^{(v-j)sp(t)/\sigma }\Big(%
\int_{2^{j}}^{2^{j+1}}\eta _{\tau }\frac{d\tau }{\tau }\Big)^{p(t)}\Big)%
^{1/\sigma },
\end{equation*}%
where $c>0$ is independent of $v$ and $j$. By Lemma \ref{DHHR-estimate} we
find $m>0$ such that%
\begin{eqnarray*}
&&\Big(\frac{1}{(j-v+3)\log 2}\int_{2^{v-1}}^{2^{j+2}}\eta _{\tau }\chi
_{\lbrack 2^{j},2^{1+j}]}(\tau )\frac{d\tau }{\tau }\Big)^{p(t)} \\
&\lesssim &\frac{1}{j-v+3}\int_{2^{v-1}}^{2^{j+2}}\eta _{\tau }^{p(0)}\chi
_{\lbrack 2^{j},2^{1+j}]}(\tau )\frac{d\tau }{\tau }+2^{jm}\chi _{\{t:q(t){<}%
q(0)\}}(t) \\
&\lesssim &\frac{1}{j-v+3}\int_{2^{j}}^{2^{j+2}}\eta _{\tau }^{p(0)}\frac{%
d\tau }{\tau }+2^{jm}\chi _{\{t:q(t){<}q(0)\}}(t)
\end{eqnarray*}%
for any $v-1\leq j\leq -1$ and any $t\in \lbrack 2^{v},2^{v+1}]\subset
\lbrack 2^{v-1},2^{j+2}]$. Therefore,%
\begin{equation*}
\eta _{t,1,v}^{p(t)}\lesssim \sum_{j=v-1}^{-1}2^{(v-j)sp^{-}/\sigma
}(j-v+3)^{p^{+}-1}\int_{2^{j}}^{2^{j+1}}\eta _{\tau }^{p(0)}\frac{d\tau }{%
\tau }+h_{v}
\end{equation*}%
for any $t\in \lbrack 2^{v},2^{v+1}]$, where%
\begin{equation*}
h_{v}=\sum_{j=v-1}^{-1}2^{(v-j)sp^{-}/\sigma }(j-v+3)^{p^{+}}2^{jm}.
\end{equation*}%
We have $\int_{2^{v}}^{2^{v+1}}\frac{dt}{t}\lesssim 1$. Therefore,%
\begin{equation*}
\int_{2^{v}}^{2^{1+v}}\eta _{t,1,v}^{p(t)}\frac{dt}{t}\lesssim
\sum_{j=v-1}^{-1}2^{s(v-j)p^{-}/\sigma
}(j-v+3)^{p^{+}-1}\int_{2^{j}}^{2^{j+1}}\eta _{\tau }^{p(0)}\frac{d\tau }{%
\tau }+h_{v}.
\end{equation*}%
Applying Lemma \ref{Hardy-inequality} we get%
\begin{equation*}
\sum_{v=-\infty }^{0}\int_{2^{v}}^{2^{v+1}}\eta _{t,1,v}^{p(t)}\frac{dt}{t}%
\lesssim \sum_{j=-\infty }^{-1}\int_{2^{j}}^{2^{j+1}}\eta _{\tau }^{p(0)}%
\frac{d\tau }{\tau }+c\lesssim \int_{0}^{1}\eta _{\tau }^{p(0)}\frac{d\tau }{%
\tau }+c\lesssim 1,
\end{equation*}%
by taking $m$ large enough such that $m>0$.

\textit{Estimation of }$\eta _{t,2,v}$. Let $\sigma >0$\ be such that $%
p^{+}<\sigma $.\ Again, we have%
\begin{equation*}
\Big(\sum_{j=0}^{\infty }\int_{2^{j}}^{2^{j+1}}\tau ^{-s}\eta _{\tau }\frac{%
d\tau }{\tau }\Big)^{p(t)/\sigma }\leq 2^{-vsp(t)/\sigma }\sum_{j=0}^{\infty
}2^{(v-j)sp(t)/\sigma }\Big(\int_{2^{j}}^{2^{j+1}}\eta _{\tau }\frac{d\tau }{%
\tau }\Big)^{p(t)/\sigma }.
\end{equation*}%
By H\"{o}lder's inequality, we estimate this expression by%
\begin{equation*}
c2^{-vsp(t)/\sigma }\Big(\sum_{j=0}^{\infty }2^{(v-j)sp(t)/\sigma }\Big(%
\int_{2^{j}}^{2^{j+1}}\eta _{\tau }\frac{d\tau }{\tau }\Big)^{p(t)}\Big)%
^{1/\sigma }.
\end{equation*}%
Again, by Lemma \ref{DHHR-estimate} we find $m>0$ such that%
\begin{eqnarray*}
&&\Big(\frac{1}{(j-v+1)\log 2}\int_{2^{v}}^{2^{j+1}}\eta _{\tau }\chi
_{\lbrack 2^{j},2^{1+j}]}(\tau )\frac{d\tau }{\tau }\Big)^{p(t)} \\
&\lesssim &\frac{1}{j-v+1}\int_{2^{v}}^{2^{j+1}}\eta _{\tau }^{p_{\infty
}}\chi _{\lbrack 2^{j},2^{1+j}]}(\tau )\frac{d\tau }{\tau }+1 \\
&\lesssim &\frac{1}{j-v+1}\int_{2^{j}}^{2^{j+1}}\eta _{\tau }^{p_{\infty }}%
\frac{d\tau }{\tau }+1
\end{eqnarray*}%
for any $j\geq 0$ and any $t\in \lbrack 2^{v},2^{v+1}]\subset \lbrack
2^{v},2^{j+1}]$. Therefore,%
\begin{equation*}
\eta _{t,2,v}^{p(t)}\lesssim \sum_{j=0}^{\infty }2^{(v-j)sp^{-}/\sigma
}(j-v+1)^{p^{+}-1}\int_{2^{j}}^{2^{j+1}}\eta _{\tau }^{p_{\infty }}\frac{%
d\tau }{\tau }+h_{v}
\end{equation*}%
for any $t\in \lbrack 2^{v},2^{v+1}]$, where%
\begin{equation*}
h_{v}=\sum_{j=0}^{\infty }2^{(v-j)sp^{-}/\sigma }(j-v+1)^{p^{+}}.
\end{equation*}%
We have $\int_{2^{v}}^{2^{v+1}}\frac{dt}{t}\lesssim 1$. Observe that%
\begin{equation*}
h_{v}\leq 2^{\frac{vsp^{-}}{2\sigma }}\sum_{j=0}^{\infty
}2^{(v-j)sp^{-}/2\sigma }(j-v+1)^{p^{+}}\lesssim 2^{\frac{vsp^{-}}{2\sigma }%
},\quad v\leq 0.
\end{equation*}%
Therefore,%
\begin{equation*}
\int_{2^{v}}^{2^{v+1}}\eta _{t,2,v}^{p(t)}\frac{dt}{t}\lesssim
\sum_{j=0}^{\infty }2^{(v-j)sp^{-}/\sigma
}(j-v+1)^{p^{+}-1}\int_{2^{j}}^{2^{j+1}}\eta _{\tau }^{p_{\infty }}\frac{%
d\tau }{\tau }+2^{\frac{v}{2\sigma }}.
\end{equation*}%
Again, by Lemma \ref{Hardy-inequality} we get%
\begin{equation*}
\sum_{v=-\infty }^{-1}\int_{2^{v}}^{2^{v+1}}\eta _{t,2,v}^{p(t)}\frac{dt}{t}%
\lesssim \sum_{j=0}^{\infty }\int_{2^{j}}^{2^{j+1}}\eta _{\tau }^{p_{\infty
}}\frac{d\tau }{\tau }+c\lesssim \int_{1}^{\infty }\eta _{\tau }^{p_{\infty
}}\frac{d\tau }{\tau }+c\lesssim 1.
\end{equation*}

\textit{Substep 1.2.} $v>0$. Let $\sigma >0$\ be such that $p^{+}<\sigma $.\
We have%
\begin{equation*}
\Big(\sum_{j=v-1}^{\infty }\int_{2^{j}}^{2^{j+1}}\tau ^{-s}\eta _{\tau }%
\frac{d\tau }{\tau }\Big)^{p(t)/\sigma }\leq 2^{-vsp(t)/\sigma
}\sum_{j=v-1}^{\infty }2^{(v-j)sp(t)/\sigma }\Big(\int_{2^{j}}^{2^{j+1}}\eta
_{\tau }\frac{d\tau }{\tau }\Big)^{p(t)/\sigma }.
\end{equation*}%
By H\"{o}lder's inequality, we estimate this expression by%
\begin{equation*}
c2^{-vsp(t)/\sigma }\Big(\sum_{j=v-1}^{\infty }2^{(v-j)sp(t)/\sigma }\Big(%
\int_{2^{j}}^{2^{j+1}}\eta _{\tau }\frac{d\tau }{\tau }\Big)^{p(t)}\Big)%
^{1/\sigma }.
\end{equation*}%
Applying Lemma \ref{DHHR-estimate} we find $m>0$ such that%
\begin{eqnarray*}
&&\Big(\frac{1}{(j-v+3)\log 2}\int_{2^{v-1}}^{2^{j+2}}\eta _{\tau }\chi
_{\lbrack 2^{j},2^{1+j}]}(\tau )\frac{d\tau }{\tau }\Big)^{p(t)} \\
&\lesssim &\frac{1}{j-v+3}\int_{2^{v-1}}^{2^{j+2}}\eta _{\tau }^{p_{\infty
}}\chi _{\lbrack 2^{j},2^{1+j}]}(\tau )\frac{d\tau }{\tau }+2^{-vm} \\
&\lesssim &\frac{1}{j-v+3}\int_{2^{j}}^{2^{j+1}}\eta _{\tau }^{p_{\infty }}%
\frac{d\tau }{\tau }+2^{-vm}
\end{eqnarray*}%
for any $j\geq v-1$ and any $t\in \lbrack 2^{v},2^{v+1}]\subset \lbrack
2^{v-1},2^{j+2}]$. Therefore,%
\begin{equation*}
\eta _{t,v}^{p(t)}\lesssim \sum_{j=v-1}^{\infty }2^{(v-j)sp^{-}/\sigma
}(j-v+2)^{p^{+}-1}\int_{2^{j}}^{2^{j+1}}\eta _{\tau }^{p_{\infty }}\frac{%
d\tau }{\tau }+h_{v}
\end{equation*}%
for any $t\in \lbrack 2^{v},2^{v+1}]$, where%
\begin{equation*}
h_{v}=\sum_{j=v-1}^{\infty }2^{(v-j)sp^{-}/\sigma }(j-v+3)^{p^{+}}2^{-vm}.
\end{equation*}%
Therefore,%
\begin{equation*}
\int_{2^{v}}^{2^{v+1}}\eta _{t,v}^{p(t)}\frac{dt}{t}\lesssim
\sum_{j=v-1}^{\infty
}2^{(v-j)sp^{-}}(j-v+3)^{p^{+}-1}\int_{2^{j}}^{2^{j+1}}\eta _{\tau
}^{p_{\infty }}\frac{d\tau }{\tau }+h_{v}.
\end{equation*}%
Applying Lemma \ref{Hardy-inequality} we get%
\begin{equation*}
\sum_{v=1}^{\infty }\int_{2^{v}}^{2^{v+1}}\eta _{t,v}^{p(t)}\frac{dt}{t}%
\lesssim \sum_{j=0}^{\infty }\int_{2^{j}}^{2^{1+j}}\eta _{\tau }^{p_{\infty
}}\frac{d\tau }{\tau }+c\lesssim \int_{1}^{\infty }\eta _{\tau }^{p_{\infty
}}\frac{d\tau }{\tau }+c\lesssim 1,
\end{equation*}%
by taking $m$ large enough such that $m>0$. The proof is completed by the
scaling argument.

\textit{Step 2.} We prove that%
\begin{equation}
\Big(\int_{0}^{1}\eta _{t}^{p(0)}\frac{dt}{t}\Big)^{\frac{1}{p(0)}}+\Big(%
\int_{1}^{\infty }\eta _{t}^{p_{\infty }}\frac{dt}{t}\Big)^{\frac{1}{%
p_{\infty }}}\lesssim \left\Vert \eta _{t}\right\Vert _{L^{p(\cdot
)}((0,\infty ),\frac{dt}{t})}.  \label{Est1.1}
\end{equation}%
We suppose that the right-hand side is less than or equal one. We will prove
that 
\begin{equation*}
\sum_{v=1}^{\infty }\int_{2^{-v}}^{2^{1-v}}\eta _{t}^{p(0)}\frac{dt}{t}%
\lesssim 1.
\end{equation*}%
Clearly follows from the inequality%
\begin{equation*}
\eta _{t}^{p(0)}\lesssim \int_{2^{-v-1}}^{2^{1-v}}\eta _{\tau }^{p(\tau )}%
\frac{d\tau }{\tau }+2^{-v}=\delta
\end{equation*}%
for any $v\in \mathbb{N}$ and any any $t\in \lbrack 2^{-v},2^{1-v}]$. This
claim can be reformulated as showing that%
\begin{equation*}
\big(\delta ^{-\frac{1}{p(0)}}\eta _{t}\big)^{p(0)}\leq \Big(\frac{1}{\log 2}%
\int_{2^{-v-1}}^{2^{1-v}}\delta ^{-\frac{1}{p(0)}}\eta _{\tau }\frac{d\tau }{%
\tau }\Big)^{p(0)}\lesssim 1.
\end{equation*}%
By Lemma \ref{DHHR-estimate},%
\begin{equation*}
\Big(\frac{\gamma }{\log 2}\int_{2^{-v-1}}^{2^{1-v}}\delta ^{-\frac{1}{p(0)}%
}\eta _{\tau }\frac{d\tau }{\tau }\Big)^{p\left( t\right) }\lesssim
\int_{2^{-v-1}}^{2^{1-v}}\delta ^{-\frac{p\left( \tau \right) }{p(0)}}\eta
_{\tau }^{p\left( \tau \right) }\frac{d\tau }{\tau }+1,
\end{equation*}%
where $\gamma =e^{-4mc_{\log }(1/p)}$ and $m>0$. We use the log-H\"{o}lder
continuity of $p$ at the origin to show that%
\begin{equation*}
\delta ^{-\frac{p\left( \tau \right) }{p(0)}}\approx \delta ^{-1},\text{ \ \ 
}\tau \in \lbrack 2^{-v-1},2^{1-v}],v\in \mathbb{N}.
\end{equation*}%
Therefore, from the definition of $\delta $, we find that%
\begin{equation*}
\int_{2^{-v-1}}^{2^{1-v}}\delta ^{-1}\eta _{\tau }^{p\left( \tau \right) }%
\frac{d\tau }{\tau }\lesssim 1
\end{equation*}%
for any $v\in \mathbb{N}$ and this implies that%
\begin{equation*}
\big(\delta ^{-\frac{1}{p(0)}}\eta _{t}\big)^{p(0)}\lesssim 1
\end{equation*}%
for any $v\in \mathbb{N}$ and any any $t\in \lbrack 2^{-v},2^{1-v}]$.

Now, we will prove that 
\begin{equation*}
\sum_{v=1}^{\infty }\int_{2^{v}}^{2^{v+1}}\eta _{t}^{p_{\infty }}\frac{dt}{t}%
\lesssim 1.
\end{equation*}%
Clearly follows from the inequality%
\begin{equation*}
\eta _{t}^{p_{\infty }}\lesssim \int_{2^{v-1}}^{2^{v+1}}\eta _{\tau
}^{p(\tau )}\frac{d\tau }{\tau }+2^{-v}=\delta
\end{equation*}%
for any $v\in \mathbb{N}$ and any any $t\in \lbrack 2^{v},2^{v+1}]$. This
claim can be reformulated as showing that%
\begin{equation*}
\big(\delta ^{-\frac{1}{p_{\infty }}}\eta _{t}\big)^{p_{\infty }}\leq \Big(%
\frac{1}{\log 2}\int_{2^{v-1}}^{2^{v+1}}\delta ^{-\frac{1}{p_{\infty }}}\eta
_{\tau }\frac{d\tau }{\tau }\Big)^{p_{\infty }}\lesssim 1.
\end{equation*}%
By Lemma \ref{DHHR-estimate},%
\begin{equation*}
\Big(\frac{\gamma }{\log 2}\int_{2^{v-1}}^{2^{v+1}}\delta ^{-\frac{1}{%
p_{\infty }}}\eta _{\tau }\frac{d\tau }{\tau }\Big)^{p\left( t\right)
}\lesssim \int_{2^{v-1}}^{2^{v+1}}\delta ^{-\frac{p\left( \tau \right) }{%
p_{\infty }}}\eta _{\tau }^{p\left( \tau \right) }\frac{d\tau }{\tau }+1,
\end{equation*}%
where $\gamma =e^{-4mc_{\log }}$ and $m>0$. We use the logarithmic decay
condition on $q$ at infinity to show that%
\begin{equation*}
\delta ^{-\frac{p\left( \tau \right) }{p_{\infty }}}\approx \delta ^{-1},%
\text{ \ \ }\tau \in \lbrack 2^{v-1},2^{v+1}],v\in \mathbb{N}.
\end{equation*}%
Therefore, from the definition of $\delta $, we find that%
\begin{equation*}
\int_{2^{v-1}}^{2^{v+1}}\delta ^{-1}\eta _{\tau }^{p\left( \tau \right) }%
\frac{d\tau }{\tau }\lesssim 1
\end{equation*}%
for any $v\in \mathbb{N}$ and this implies that for any $v\in \mathbb{N}$
any any $t\in \lbrack 2^{v},2^{v+1}]$, 
\begin{equation*}
\big(\delta ^{-\frac{1}{p_{\infty }}}\eta _{t}\big)^{p_{\infty }}\lesssim 1,
\end{equation*}%
which\ completes\ the proof of $\mathrm{\eqref{Est1.1}}$, by the scaling
argument.

\textit{Step 3.} We prove that%
\begin{equation}
\left\Vert \lambda _{t}\right\Vert _{L^{p(\cdot )}((0,\infty ),\frac{dt}{t}%
)}\lesssim \Big(\int_{0}^{1}\lambda _{t}^{p(0)}\frac{dt}{t}\Big)^{\frac{1}{%
p(0)}}+\Big(\int_{1}^{\infty }\lambda _{t}^{p_{\infty }}\frac{dt}{t}\Big)^{%
\frac{1}{p_{\infty }}}.  \label{Est2}
\end{equation}%
We suppose that the right-hand side is less than or equal one. Notice that%
\begin{equation*}
\left\Vert \lambda _{t}\right\Vert _{L^{p(\cdot )}((0,\infty ),\frac{dt}{t}%
)}\approx \Big\|\Big(t^{-\frac{1}{p(t)}}\lambda _{t}\chi _{\lbrack
2^{v},2^{1+v}]}\Big)_{v}\Big\|_{\ell ^{p(\cdot )}(L^{p(\cdot )})}.
\end{equation*}%
We see that%
\begin{eqnarray*}
\lambda _{t} &=&\frac{1}{\log 2}\int_{t}^{2t}\lambda _{t}\frac{d\tau }{\tau }%
\leq \frac{1}{\log 2}\int_{t}^{2t}\lambda _{\tau }\frac{d\tau }{\tau }%
\lesssim t^{-s}\int_{0}^{2t}\tau ^{s}\lambda _{\tau }\frac{d\tau }{\tau } \\
&\leq &t^{-s}\int_{0}^{2^{v+2}}\tau ^{s}\lambda _{\tau }\frac{d\tau }{\tau }
\\
&\leq &\sum_{i=-\infty }^{v}t^{-s}\int_{2^{i}}^{2^{i+2}}\tau ^{s}\lambda
_{\tau }\frac{d\tau }{\tau }=\sum_{j=-v}^{\infty
}t^{-s}\int_{2^{-j}}^{2^{2-j}}\tau ^{s}\lambda _{\tau }\frac{d\tau }{\tau }
\end{eqnarray*}%
for any $v\leq 0$ any any $t\in \lbrack 2^{v},2^{v+1}]$. Let $\sigma >0$ be
such that $p^{+}<\sigma $. We have%
\begin{equation*}
\Big(\sum_{j=-v}^{\infty }\int_{2^{-j}}^{2^{2-j}}\tau ^{s}\lambda _{\tau }%
\frac{d\tau }{\tau }\Big)^{p(t)/\sigma }\leq 2^{vsp(t)/\sigma
}\sum_{j=-v}^{\infty }2^{-(v+j)sp(t)/\sigma }\Big(\int_{2^{-j}}^{2^{2-j}}%
\lambda _{\tau }\frac{d\tau }{\tau }\Big)^{p(t)/\sigma }.
\end{equation*}%
Again, by H\"{o}lder's inequality, we estimate this expression by%
\begin{equation*}
c2^{vsp(t)/\sigma }\Big(\sum_{j=-v}^{\infty }2^{-(v+j)sp(t)/\sigma }\Big(%
\int_{2^{-j}}^{2^{2-j}}\lambda _{\tau }\frac{d\tau }{\tau }\Big)^{p(t)}\Big)%
^{1/\sigma }.
\end{equation*}%
Applying again Lemma \ref{DHHR-estimate} we get%
\begin{eqnarray*}
&&\Big(\frac{1}{(j+v+2)\log 2}\int_{2^{-j}}^{2^{v+2}}\lambda _{\tau }\chi
_{\lbrack 2^{-j},2^{2-j}]}(\tau )\frac{d\tau }{\tau }\Big)^{p(t)} \\
&\lesssim &\frac{1}{j+v+2}\int_{2^{-j}}^{2^{v+2}}\lambda _{\tau }^{p(0)}\chi
_{\lbrack 2^{-j},2^{2-j}]}(\tau )\frac{d\tau }{\tau }+2^{vm}.
\end{eqnarray*}%
Therefore,%
\begin{equation*}
\lambda _{t}^{p(t)}\lesssim \sum_{j=-v}^{\infty }2^{-(v+j)sp(t)/\sigma
}(j+v+2)^{p^{+}-1}\int_{2^{-j}}^{2^{2-j}}\lambda _{\tau }^{p(0)}\frac{d\tau 
}{\tau }+f_{v}
\end{equation*}%
for $v\leq 0$ and any $t\in \lbrack 2^{v},2^{v+1}]\subset \lbrack
2^{-j},2^{v+2}]$\ where%
\begin{equation*}
f_{v}=2^{vm}.
\end{equation*}%
Therefore,%
\begin{equation*}
\int_{2^{v}}^{2^{v+1}}\lambda _{t}^{p(t)}\frac{dt}{t}\lesssim
\sum_{j=-v}^{\infty }2^{-(v+j)sp^{-}/\sigma
}(j+v+2)^{p^{+}-1}\int_{2^{-j}}^{2^{2-j}}\lambda _{\tau }^{p(0)}\frac{d\tau 
}{\tau }+f_{v}.
\end{equation*}%
By taking $m$ large enough such that $m>0$ and again\ by Lemma \ref%
{Hardy-inequality}\ we\ get%
\begin{equation*}
\sum_{v=-\infty }^{0}\int_{2^{v}}^{2^{v+1}}\lambda _{t}^{p(t)}\frac{dt}{t}%
\lesssim \sum_{j=1}^{\infty }\int_{2^{-j}}^{2^{2-j}}\lambda _{\tau }^{p(0)}%
\frac{d\tau }{\tau }+c\lesssim 1.
\end{equation*}%
Now we see that%
\begin{eqnarray*}
\lambda _{t} &=&\int_{t}^{2t}\lambda _{t}\frac{d\tau }{\tau }\leq
\int_{t}^{2t}\lambda _{\tau }\frac{d\tau }{\tau }\lesssim
t^{-s}\int_{1}^{2t}\tau ^{s}\lambda _{\tau }\frac{d\tau }{\tau } \\
&\leq &t^{-s}\int_{1}^{2^{2+v}}\tau ^{s}\lambda _{\tau }\frac{d\tau }{\tau }
\\
&\leq &\sum_{j=0}^{v}t^{-s}\int_{2^{j}}^{2^{j+2}}\tau ^{s}\lambda _{\tau }%
\frac{d\tau }{\tau }
\end{eqnarray*}%
for any $v>0$ any any $t\in \lbrack 2^{v},2^{v+1}]$. Let $\sigma >0$ be such
that $p^{+}<\sigma $.\ We have%
\begin{eqnarray*}
\Big(\sum_{j=0}^{v}\int_{2^{j}}^{2^{j+2}}\tau ^{s}\lambda _{\tau }\frac{%
d\tau }{\tau }\Big)^{p(t)/\sigma } &\leq &\sum_{j=0}^{v}\Big(%
\int_{2^{j}}^{2^{j+2}}\tau ^{s}\lambda _{\tau }\frac{d\tau }{\tau }\Big)%
^{p(t)/\sigma } \\
&\leq &\sum_{j=0}^{v}2^{jsp(t)/\sigma }\Big(\int_{2^{j}}^{2^{j+2}}\lambda
_{\tau }\frac{d\tau }{\tau }\Big)^{p(t)/\sigma } \\
&=&2^{vsp(t)/\sigma }\sum_{j=0}^{v}2^{(j-v)sp(t)/\sigma }\Big(%
\int_{2^{j}}^{2^{j+2}}\lambda _{\tau }\frac{d\tau }{\tau }\Big)^{p(t)/\sigma
}.
\end{eqnarray*}%
Again, by H\"{o}lder's inequality, we estimate this expression by%
\begin{equation*}
c2^{vsp(t)/\sigma }\Big(\sum_{j=0}^{v}2^{(j-v)sp(t)/\sigma }\Big(%
\int_{2^{j}}^{2^{j+2}}\lambda _{\tau }\frac{d\tau }{\tau }\Big)^{p(t)}\Big)%
^{1/\sigma }.
\end{equation*}%
Applying again Lemma \ref{DHHR-estimate} we get%
\begin{eqnarray*}
&&\Big(\frac{1}{(v-j+2)\log 2}\int_{2^{j}}^{2^{v+2}}\lambda _{\tau }\chi
_{\lbrack 2^{j},2^{2+j}]}(\tau )\frac{d\tau }{\tau }\Big)^{p(t)} \\
&\lesssim &\frac{1}{v-j+2}\int_{2^{j}}^{2^{v+2}}\lambda _{\tau }^{p_{\infty
}}\chi _{\lbrack 2^{j},2^{2+j}]}(\tau )\frac{d\tau }{\tau }+2^{-jm}
\end{eqnarray*}%
for $v>0$ and any $t\in \lbrack 2^{v},2^{v+1}]\subset \lbrack 2^{j},2^{v+1}]$%
. Therefore,%
\begin{equation*}
\lambda _{t}^{p(t)}\lesssim \sum_{j=0}^{v}2^{(j-v)sp(t)/\sigma
}(v-j+2)^{p^{+}-1}\int_{2^{j}}^{2^{j+2}}\lambda _{\tau }^{p_{\infty }}\frac{%
d\tau }{\tau }+f_{v}
\end{equation*}%
for $v>0$ and any $t\in \lbrack 2^{v},2^{v+1}]$, where%
\begin{equation*}
f_{v}=\sum_{j=0}^{v}2^{(j-v)sp^{-}/\sigma }(v-j+2)^{p^{+}}2^{-jm}.
\end{equation*}%
Therefore,%
\begin{equation*}
\int_{2^{v}}^{2^{v+1}}\lambda _{t}^{p(t)}\frac{dt}{t}\lesssim
\sum_{j=0}^{v}2^{(j-v)sp^{-}/\sigma
}(v-j+2)^{p^{+}-1}\int_{2^{j}}^{2^{j+2}}\lambda _{\tau }^{p_{\infty }}\frac{%
d\tau }{\tau }+f_{v}.
\end{equation*}%
By taking $m$ large enough such that $m>0$ and again\ by Lemma \ref%
{Hardy-inequality}\ we\ get%
\begin{equation*}
\sum_{v=1}^{\infty }\int_{2^{v}}^{2^{v+1}}\lambda _{t}^{p(t)}\frac{dt}{t}%
\lesssim \sum_{j=0}^{\infty }\int_{2^{j}}^{2^{2+j}}\lambda _{\tau
}^{p_{\infty }}\frac{d\tau }{\tau }+c\lesssim 1.
\end{equation*}%
The proof of $\mathrm{\eqref{Est2}}$ is completed by the scaling argument.

\textit{Step 4.} We prove that%
\begin{equation*}
\Big(\int_{0}^{1}\lambda _{t}^{p(0)}\frac{dt}{t}\Big)^{\frac{1}{p(0)}}+\Big(%
\int_{1}^{\infty }\lambda _{t}^{p_{\infty }}\frac{dt}{t}\Big)^{\frac{1}{%
p_{\infty }}}\lesssim \left\Vert \lambda _{t}\right\Vert _{L^{p(\cdot
)}((0,\infty ),\frac{dt}{t})}.
\end{equation*}

We omit the proofs of this estimate, since they are essentially similar to
the proof of $\mathrm{\eqref{Est1.1}}$. \quad $\square $

We would like to mention that the estimates $\mathrm{\eqref{Est1.1.1}}$ and $%
\mathrm{\eqref{Est1.1.2}}$ are true if we assume that%
\begin{equation*}
\eta _{t}\leq \eta _{\tau },\quad 0<\tau \leq t
\end{equation*}%
and%
\begin{equation*}
\lambda _{t}\leq \lambda _{\tau },\quad 0<t\leq \tau ,
\end{equation*}%
respectively. Also we find that 
\begin{equation*}
\left\Vert \eta _{t}\right\Vert _{L^{q(\cdot )}((0,\infty ),\frac{dt}{t}%
)}\approx \left\Vert \eta _{t}\right\Vert _{L^{p(\cdot )}((0,\infty ),\frac{%
dt}{t})}
\end{equation*}%
and%
\begin{equation*}
\left\Vert \lambda _{t}\right\Vert _{L^{q(\cdot )}((0,\infty ),\frac{dt}{t}%
)}\approx \left\Vert \lambda _{t}\right\Vert _{L^{p(\cdot )}((0,\infty ),%
\frac{dt}{t})}
\end{equation*}%
for any $p,q\in \mathcal{P}(\mathbb{R})$\textit{\ }are \emph{$\log $}-H\"{o}%
lder continuous both at the origin\ and at infinity with $1\leq q^{-}\leq
q^{+}<\infty $, $1\leq p^{-}\leq p^{+}<\infty $, 
\begin{equation*}
p(0)=q(0)\quad \text{and\quad }p_{\infty }=q_{\infty }.
\end{equation*}

By the technical of Theorem \ref{main-result} we immediately arrive at the
following result.

\begin{theorem}
\label{main-result1}\textit{Let }$s>0$\textit{. Let }$p\in \mathcal{P}(%
\mathbb{R})$\textit{\ }be \emph{$\log $}-H\"{o}lder continuous both at the
origin with $1\leq p^{-}\leq p^{+}<\infty $\textit{. Let }$\left\{
\varepsilon _{t}\right\} _{t}$\textit{\ be a sequence of positive measurable
functions.} Let%
\begin{equation*}
\eta _{t}=t^{s}\int_{t}^{1}\tau ^{-s}\varepsilon _{\tau }\frac{d\tau }{\tau }%
\quad \text{and\quad }\lambda _{t}=t^{-s}\int_{0}^{t}\tau ^{s}\varepsilon
_{\tau }\frac{d\tau }{\tau }.
\end{equation*}%
Then there exists constant $c>0\ $\textit{depending only on }$s$, $p^{-}$%
\textit{, c}$_{\log }(p)$ \textit{and }$p^{+}$ such that%
\begin{equation}
\left\Vert \eta _{t}\right\Vert _{L^{p(\cdot )}((0,1],\frac{dt}{t})}\approx %
\Big(\int_{0}^{1}\eta _{t}^{p(0)}\frac{dt}{t}\Big)^{\frac{1}{p(0)}}
\label{Est1.1.3}
\end{equation}%
and%
\begin{equation}
\left\Vert \lambda _{t}\right\Vert _{L^{p(\cdot )}((0,1],\frac{dt}{t}%
)}\approx \Big(\int_{0}^{1}\lambda _{t}^{p(0)}\frac{dt}{t}\Big)^{\frac{1}{%
p(0)}}.  \label{Est1.1.4}
\end{equation}%
Moreover, 
\begin{equation*}
\left\Vert \eta _{t}\right\Vert _{L^{p(\cdot )}((0,1],\frac{dt}{t}%
)}+\left\Vert \lambda _{t}\right\Vert _{L^{p(\cdot )}((0,1],\frac{dt}{t}%
)}\lesssim \Big(\int_{0}^{1}\varepsilon _{t}^{p(0)}\frac{dt}{t}\Big)^{\frac{1%
}{p(0)}}.
\end{equation*}
\end{theorem}

Again, we would like to mention that the estimates $\mathrm{\eqref{Est1.1.3}}
$ and $\mathrm{\eqref{Est1.1.4}}$ are true if we assume that%
\begin{equation*}
\eta _{t}\leq \eta _{\tau },\quad 0<\tau \leq t\leq 1
\end{equation*}%
and%
\begin{equation*}
\lambda _{t}\leq \lambda _{\tau },\quad 0<t\leq \tau \leq 1,
\end{equation*}%
respectively. Also we find that 
\begin{equation*}
\left\Vert \eta _{t}\right\Vert _{L^{q(\cdot )}((0,1],\frac{dt}{t})}\approx
\left\Vert \eta _{t}\right\Vert _{L^{p(\cdot )}((0,1],\frac{dt}{t})}
\end{equation*}%
and%
\begin{equation*}
\left\Vert \lambda _{t}\right\Vert _{L^{q(\cdot )}((0,1],\frac{dt}{t}%
)}\approx \left\Vert \lambda _{t}\right\Vert _{L^{p(\cdot )}((0,1],\frac{dt}{%
t})}
\end{equation*}%
for any $p,q\in \mathcal{P}(\mathbb{R})$\textit{\ }are \emph{$\log $}-H\"{o}%
lder continuous at the origin with $1\leq q^{-}\leq q^{+}<\infty $, $1\leq
p^{-}\leq p^{+}<\infty $ and 
\begin{equation*}
p(0)=q(0).
\end{equation*}

Douadi Drihem

M'sila University, Department of Mathematics,

Laboratory of Functional Analysis and Geometry of Spaces,

P.O. Box 166, M'sila 28000, Algeria,

e-mail: \texttt{\ douadidr@yahoo.fr}

\end{document}